\documentclass[12pt]{amsart}
\newtheorem{thm}{Theorem}
\newtheorem{lemma}{Lemma}
\newtheorem{cor}{Corollary}

\newcommand{\rk}{\operatorname{rank}}
\newcommand{\ffrac}[2]{{\mbox{\large$\frac{#1}{#2}$}}}

\begin{document}
\title{CR geometry and conformal foliations}
\author{Paul Baird}
\address{\hskip-\parindent
D\'epartement de Math\'ematiques\\
Universit\'e de Bretagne Occidentale\newline
6 av.\ Victor Le Gorgeu -- CS 93837\\
29238 Brest Cedex\\
France}
\email{Paul.Baird@univ-brest.fr}
\author{Michael Eastwood}
\address{\hskip-\parindent
Mathematical Sciences Institute\\
Australian National University,\newline
ACT 0200, Australia}
\email{meastwoo@member.ams.org}

\thanks{The first author is grateful to the Australian National University for
support and hospitality while working on this article. The second author is a
Federation Fellow of the Australian Research Council.}
\begin{abstract}
We use the CR geometry of the standard hyperquadric in ${\mathbb{CP}}_3$ to
give a detailed twistor description of conformal foliations in Euclidean
3-space.
\end{abstract}
\renewcommand{\subjclassname}{\textup{2010} Mathematics Subject Classification}
\maketitle

\section{Introduction}
A foliation of a Riemannian manifold is said to be conformal if and only if the
corresponding locally defined submersion is conformal on the orthogonal spaces
to the leaves (precise definitions are given in~\S\ref{conformalfoliations}).
Locally a conformal foliation induces
${\mathbb{R}}^3\supseteq{}^{\mathrm{open}}\Omega\stackrel{h}{\to}{\mathbb{C}}$,
defined up to composition $\alpha\circ h$ for $\alpha$ conformal, with the
property that the gradient of $h$ is null: $(\nabla h)^2=0$. Such a mapping is
said to be horizontally conformal~\cite{thebook}. In~\cite{n} Nurowski showed
how real-analytic horizontally conformal mappings in ${\mathbb{R}}^3$ may be
constructed from a holomorphic function of two complex variables. The mappings
obtained in this way are certainly real-analytic. Otherwise they are generic 
but, as we shall see, not completely general.

In this article we construct conformal foliations themselves starting with a
holomorphic function of two complex variables. The foliations constructed in
this way are real-analytic but are otherwise general, at least locally. Our
construction is one of twistor geometry, a viewpoint which also allows us to
identify the smooth conformal foliations with local CR hypersurfaces in the
standard Levi-indefinite hyperquadric in~${\mathbb{CP}}_3$. Real-analyticity is
a familiar distinguishing feature in CR geometry. Our construction parallels
the twistor interpretation of the Kerr Theorem in special relativity as
explained in \cite{NT} and, following~\cite{bw}, the precise link is explained
at the end of our paper in~\S\ref{kerr}. Another intimately related
construction appears in~\cite{bp} where twistor geometry is used directly to
describe horizontally conformal mappings in~${\mathbb{R}}^3$.

We would like to thank John Bland for several useful conversations on CR
geometry and for crucial observations concerning Lagrangian geometry as it
appears in \S\ref{constructandcompare} and especially Theorem~\ref{lagrangian}.

\section{CR geometry}\label{CRgeometry}
The book \cite{ber} by Baouendi, Ebenfelt, and Rothschild provides a good
reference for CR geometry. Here, we provide an outline of the specific results
we shall need, referring to~\cite{ber} for proofs and further detail.

We start with some linear algebra, presented in the dimensions where it will be
needed. Let us consider a real linear subspace $T\subseteq{\mathbb{C}}^3$ of
real codimension~$d$. Let us denote by $J$ the real linear endomorphism of
${\mathbb{C}}^3$ given by multiplication by $i$ and let $H=T\cap JT$. It is 
the maximal complex-linear subspace of~$T$.
\begin{itemize}
\item Case $d=1$: it follows that $\dim_{\mathbb{C}}H=2$ and all such $T$ are 
on an equal footing. More precisely, ${\mathrm{GL}}(3,{\mathbb{C}})$ acts 
transitively on ${\mathrm{Gr}}_5({\mathbb{R}}^6)$, the Grassmannian of real
hyperplanes in~${\mathbb{R}}^6$. 
\item Case $d=2$: there are two cases according to whether $T$ is a complex
subspace. If not, then $\dim_{\mathbb{C}}H=1$ and we shall refer to $T$ as 
{\em generic\/}. There are two orbits for the action of
${\mathrm{GL}}(3,{\mathbb{C}})$ on ${\mathrm{Gr}}_4({\mathbb{R}}^6)$, namely
${\mathrm{Gr}}_2({\mathbb{C}}^3)$ and its complement. 
\item Case $d=3$: again there are two orbits for the action of
${\mathrm{GL}}(3,{\mathbb{C}})$ on the relevant Grassmannian. Generically 
$H=0$ and we shall refer to $T$ as {\em totally real\/}. Otherwise 
$\dim_{\mathbb{C}}H=1$. 
\end{itemize}

Now suppose $Z$ is a complex manifold with $\dim_{\mathbb{C}}Z=3$ and 
$M\subset Z$ is smooth real submanifold of real codimension~$d$. When $d=1$, 
we may apply the construction above in each tangent space, obtaining a smooth 
subbundle $H\subset TM$ equipped with an endomorphism $J:H\to H$ with 
$J^2=-{\mathrm{Id}}$. Additionally, the bundle 
\begin{equation}\label{H01}
H^{0,1}M=\{X\in {\mathbb{C}}H\mbox{ s.t.\ }JX=-iX\}\end{equation}
is closed under Lie bracket. In abstraction, such a structure is called a 
{\em CR~structure of hypersurface type\/}.

When $d=2$, we shall say that $M$ is a {\em CR~submanifold\/} if and only if
$T_xM\cap JT_xM$ is of constant dimension as $x\in M$ varies. Again, we obtain
a subbundle $H\subset TM$ equipped with a complex structure $J:H\to H$. When
$\rk_{\mathbb{C}}H=1$ we shall say that $M\subset Z$ is {\em generic\/}. If
$\rk_{\mathbb{C}}H=2$, then $M$ is simply a complex submanifold of~$Z$.

Finally, when $d=3$ there are two cases corresponding infinitesimally to the
two orbits of ${\mathrm{GL}}(3,{\mathbb{C}})$ on
${\mathrm{Gr}}_3({\mathbb{R}}^6)$ identified above. When $H=0$ we shall say
that $M$ is {\em totally real\/}. Otherwise $\rk_{\mathbb{C}}H=1$. In general,
$\rk_{\mathbb{C}}H$ is called the {\em CR dimension\/} of~$M$. 

An abstract CR structure on a smooth manifold $M$ of CR dimension $m$ is
defined by a smooth subbundle $H\subset TM$ of real rank $2m$ equipped with an 
endomorphism $J$ with $J^2=-{\mathrm{Id}}$ and such that the bundle $H^{0,1}$ 
defined by (\ref{H01}) is closed under Lie bracket. If the CR structure on $M$ 
is induced by an embedding $M\hookrightarrow Z$ in a complex manifold $Z$ and 
$f:Z\to{\mathbb{C}}$ is a holomorphic function on~$Z$, then $f|_M$ satisfies 
the partial differential equations
$$Xf=0\enskip\mbox{for all}\enskip X\in\Gamma(M,H^{0,1}).$$
These are the remnants on $M$ of the Cauchy-Riemann equations on $Z$ and
solutions of these equations on $M$ are called {\em CR functions\/} (whether or
not they arise by restriction of holomorphic functions on~$Z$).

Suppose $Q$ is a smooth CR manifold of hypersurface type and real dimension~$5$
(hence of CR dimension~$2$). Suppose that $f:Q\to{\mathbb{C}}$ is a smooth CR
function without critical points. Of course,
$$M\equiv\{p\in Q\mbox{ s.t.\ }f(p)=0\}$$
is a smooth submanifold and it is easy to check that the CR equations on $f$
imply that intrinsically $M$ is a CR manifold of CR dimension~$1$. This is the
CR analogue of the statement that holomorphic functions without critical points
on complex manifolds implicitly define complex submanifolds. Unlike the
holomorphic case, however, there is no CR implicit function theorem. As we
shall see, it is not necessarily the case that a CR submanifold of $Q$ of real
codimension 2 and CR dimension $1$ need be locally defined as the zeroes of a
CR function, even when
\begin{equation}\label{thisisQ}
Q=\{[Z_1,Z_2,Z_3,Z_4]\in{\mathbb{CP}}_3\mbox{ s.t.\ }
|Z_1|^2+|Z_2|^2=|Z_3|^2+|Z_4|^2\}\end{equation}
as we shall suppose henceforth. This particular CR manifold is known as the
hyperquadric (of CR dimension $2$ and indefinite signature). The following
theorem follows from the classical Lewy extension theorem as, for example,
proved in~\cite[Theorem~2.6.13]{h}. 

\begin{thm} Suppose $\Omega^{\mathrm{open}}\subseteq Q\subset{\mathbb{CP}}_3$
and $f:\Omega\to{\mathbb{C}}$ is a smooth CR function. Then $f$ automatically
extends as a holomorphic function to a neighbourhood of $\Omega$
in~${\mathbb{CP}}_3$. This extension is germ-unique.
\end{thm}
It is immediate that the CR functions on the indefinite hyperquadric $Q$ are
necessarily real-analytic. 
\begin{thm}\label{holomorphicextension}
Suppose $M\subset\Omega^{\mathrm{open}}\subseteq Q\subset{\mathbb{CP}}_3$ is a
real-analytic CR submanifold of real dimension $3$ and CR dimension~$1$. Then
$M$ extends into ${\mathbb{CP}}_3$ as a complex submanifold. This extension is
germ-unique. In particular, $M$ can be locally defined by a CR function.
\end{thm}
\begin{proof} Immediate from \cite[Corollary~1.8.10]{ber}. An alternative
argument may be constructed from the holomorphic implicit function theorem in
the complexification. \end{proof} 

Later in this article, we shall find smooth $3$-dimensional local CR
submanifolds of $Q$ of CR dimension 1 that are not real-analytic.

\section{Conformal foliations in ${\mathbb{R}}^3$}\label{conformalfoliations}
A unit vector field $U$ on $\Omega^{\mathrm{open}}\subseteq{\mathbb{R}}^3$
induces a $1$-dimensional foliation of $\Omega$ as its integral curves.
Conversely, all $1$-dimensional foliations arise in this way. We shall say that
$U$ is {\em transversally conformal\/} if the Lie derivative ${\mathcal{L}}_U$
preserves the conformal metric orthogonal to its integral curves. In case the
foliation is locally defined by a submersion $\pi:\Omega\to\Sigma$ then this is
equivalent to saying that $\pi$ is horizontally conformal~\cite{thebook}
(whence $\Sigma$ is naturally a Riemann surface).

By writing the Lie derivative in terms of the flat connection $\nabla$ on
${\mathbb{R}}^3$ it follows that a unit vector field $U$ on ${\mathbb{R}}^3$
defines a conformal foliation if and only if
$$\langle U,\nabla_XY+\nabla_YX\rangle=0$$
for all vector fields $X$ and $Y$ with 
$$\langle U,X\rangle=0\qquad\langle U,Y\rangle=0\qquad\langle X,Y\rangle=0
\qquad\|X\|=\|Y\|.$$ 
As detailed in~\cite{thebook}, it is equivalent to check for a particular 
non-zero pair of such fields $X,Y$, that 
$$\langle U,\nabla_XY+\nabla_YX\rangle=0\quad\mbox{and}\quad
\langle U,\nabla_XX-\nabla_YY\rangle=0.$$
We can write these conditions as partial differential equations on the
components of~$U$. Specifically, let us write $(q,r,s)$ for the standard
Euclidean co\"ordinates on~${\mathbb{R}}^3$ and write
\begin{equation}\label{thisisU}
U=\left\lgroup\!\!\begin{array}cu\\ v\\ w\end{array}\!\!\right\rgroup=
\left\lgroup\!\!\begin{array}cu(q,r,s)\\ v(q,r,s)\\ w(q,r,s)
\end{array}\!\!\right\rgroup.\end{equation}
We may suppose without loss of generality that 
$$U=\left\lgroup\!\!\begin{array}c1\\ 0\\ 0\end{array}\!\!\right\rgroup\quad
X=\left\lgroup\!\!\begin{array}c0\\ 1\\ 0\end{array}\!\!\right\rgroup\quad
Y=\left\lgroup\!\!\begin{array}c0\\ 0\\ 1\end{array}\!\!\right\rgroup\quad$$
at the origin and take 
$$X=\left\lgroup\!\!\begin{array}c-v\\ u\\ 0\end{array}\!\!\right\rgroup
\quad\mbox{and}\quad Y=
\left\lgroup\!\!\begin{array}c-uw\\ -vw\\ u^2+v^2\end{array}\!\!\right\rgroup=
\left\lgroup\!\!\begin{array}c0\\ 0\\ 1\end{array}\!\!\right\rgroup
-w\left\lgroup\!\!\begin{array}cu\\ v\\ w\end{array}\!\!\right\rgroup$$
nearby. Since $u^2+v^2+w^2=1$, it follows that all partial derivatives of $u$ 
vanish at the origin and then we compute
$$\langle U,\nabla_XY+\nabla_YX\rangle+\frac{\partial w}{\partial r}
+\frac{\partial v}{\partial s}=0=
\langle U,\nabla_XX-\nabla_YY\rangle+\frac{\partial v}{\partial r}-
\frac{\partial w}{\partial s}.$$
It follows that the conformality of $U$ is captured by the equations
\begin{equation}\label{conformality}
\frac{\partial v}{\partial r}=\frac{\partial w}{\partial s}
\quad\mbox{and}\quad
\frac{\partial v}{\partial s}=-\frac{\partial w}{\partial r}\end{equation}
at points where $(u,v,w)=(1,0,0)$.

\section{Twistor fibrations}\label{twistors} If we choose an identification 
${\mathbb{C}}^4={\mathbb{H}}^2$, where ${\mathbb{H}}$ is the space of 
quaternions, then we obtain a submersion
\begin{equation}\label{classical}
\tau:{\mathbb{CP}}_3\longrightarrow{\mathbb{HP}}_1=S^4\end{equation}
by taking the quaternionic span. See \cite{a} for details, where this fibration
is also realised as the bundle of orthogonal complex structures over the round
$4$-sphere. In co\"ordinates we have
\begin{equation}\label{formula}{\mathbb{CP}}_3\ni[Z]
\stackrel{\tau}{\mapsto}
\frac{1}{\|Z\|^2}\!\!
\left\lgroup\!\!\begin{array}c
|Z_1|^2+|Z_2|^2-|Z_3|^2-|Z_4|^2\\
2(Z_2\bar Z_3+Z_4\bar Z_1)\\
2(Z_1\bar Z_3-Z_4\bar Z_2)
\end{array}\!\!\right\rgroup\in
\begin{array}c{\mathbb{R}}\\[-2pt] \oplus\\ {\mathbb{C}}^2\end{array}
={\mathbb{R}}^5\end{equation}
and we observe that, when restricted to the hyperquadric~$Q$, this gives a 
submersion
\begin{equation}\label{cr}\tau:Q\to S^3.\end{equation}
We shall refer to both (\ref{classical}) and~(\ref{cr}) as {\em twistor 
fibrations\/}. Both have fibres intrinsically isomorphic to the Riemann sphere 
${\mathbb{CP}}_1$ with its usual complex structure.

As in~\cite{a}, it is useful to view these fibrations as intrinsically attached
to the manifolds $S^4$ and~$S^3$, regarded as flat conformal manifolds in the
usual manner. More specifically, we may use stereographic projection
$S^n\setminus\{\infty\}\stackrel{\simeq\quad}{\longrightarrow}{\mathbb{R}}^n$ to 
restrict these fibrations
$$\begin{array}{ccc}
{\mathbb{CP}}_3\setminus{\mathbb{I}}
&\stackrel{\tau}{\longrightarrow}&{\mathbb{R}}^4\\
\cup&&\cup\\
Q\setminus{\mathbb{I}}
&\stackrel{\tau}{\longrightarrow}&{\mathbb{R}}^3
\end{array}$$
where ${\mathbb{I}}=\tau^{-1}(\infty)=\{[\ast,\ast,0,0]\}$ and rewrite them as
follows.

Suppose $J\in{\mathrm{SO}}(4)$ satisfies $J^2=-{\mathrm{Id}}$. It follows that 
\begin{equation}\label{J}
J=\left\lgroup\!\!\begin{array}{cccc}
0&-u&-v&-w\\
u&0&-w&v\\
v&w&0&-u\\
w&-v&u&0
\end{array}\!\!\right\rgroup\end{equation}
for some $(u,v,w)$ s.t.\ $u^2+v^2+w^2=1$. In other words, the complex 
structures on the vector space ${\mathbb{R}}^4$ preserving the standard metric
and orientation are parameterised by $(u,v,w)\in S^2$. As in~\cite{a}, we
identify
\begin{equation}\label{cp3coords}\begin{array}{ccc}
{\mathbb{CP}}_3\setminus{\mathbb{I}}
&\stackrel{\tau}{\longrightarrow}&{\mathbb{R}}^4\\
\|&&\|\\
\makebox[0pt][r]{$(p,q,r,s,u,v,w)\in{}\,$}{\mathbb{R}}^4\times S^2
&\stackrel{\pi}{\longrightarrow}&{\mathbb{R}}^4
\end{array}\end{equation}
where $\pi$ is projection onto the first factor and we may rewrite the complex 
structure on ${\mathbb{CP}}_3$ as induced by the action of 
\begin{equation}\label{bigJ}
{\mathbb{J}}=\left\lgroup\!\!\begin{array}{ccccccc}
0&-u&-v&-w&0&0&0\\
u&0&-w&v&0&0&0\\
v&w&0&-u&0&0&0\\
w&-v&u&0&0&0&0\\
0&0&0&0&0&-w&v\\
0&0&0&0&-v&0&-u\\
0&0&0&0&-v&u&0
\end{array}\!\!\right\rgroup\end{equation}
on ${\mathbb{R}}^4\times{\mathbb{R}}^3$. Similarly, the CR structure on $Q$
when viewed as
\begin{equation}\label{Qcoords}\begin{array}{ccc}
Q\setminus{\mathbb{I}}
&\stackrel{\tau}{\longrightarrow}&{\mathbb{R}}^3\\
\|&&\|\\
\makebox[0pt][r]{$(q,r,s,u,v,w)\in{}\,$}{\mathbb{R}}^3\times S^2
&\stackrel{\pi}{\longrightarrow}&{\mathbb{R}}^3
\end{array}\end{equation}
comprises
\begin{itemize}
\item the contact structure defined by $\theta\equiv u\,dq+v\,dr+w\,ds$,
\item the endomorphism of $H\equiv\ker\theta$ induced by~${\mathbb{J}}$.
\end{itemize}

\section{Integrable Hermitian structures}
In \S\ref{twistors} we saw that the almost Hermitian structures for the
standard Euclidean metric $\Omega^{\mathrm{open}}\subseteq{\mathbb{R}}^4$ are
parameterised by smooth functions $U:{\mathbb{R}}^4\to S^2$ specifying matrices
of the form~(\ref{J}). More generally, for 
$\Omega^{\mathrm{open}}\subseteq S^4$, the almost complex structures compatible
with the flat conformal structure correspond to sections of
the twistor fibration. Let us see what it means for such complex structures to 
be integrable.
\begin{lemma}\label{nijenhuis}
At a point where $(u,v,w)=(1,0,0)$, integrability of\/
{\rm(\ref{J})} is captured by the equations
$$\frac{\partial v}{\partial p}=\frac{\partial w}{\partial q}\qquad
\frac{\partial v}{\partial q}=-\frac{\partial w}{\partial p}\qquad
\frac{\partial v}{\partial r}=\frac{\partial w}{\partial s}\qquad
\frac{\partial v}{\partial s}=-\frac{\partial w}{\partial r}$$
\end{lemma}
\begin{proof}
If we compute the Nijenhuis tensor
$$N(X,Y)\equiv[X,Y]+J[JX,Y]+J[X,JY]-[JX,JY]$$
at such a point we find that 
$$\textstyle 
N\left(\frac{\partial}{\partial p},\frac{\partial}{\partial q}\right)\!=\!
\left\lgroup\!\!\begin{array}{c}
\displaystyle
0\\ 0\\ 0\\ 0
\end{array}\!\!\right\rgroup\quad
N\left(\frac{\partial}{\partial p},\frac{\partial}{\partial r}\right)\!=\!
\left\lgroup\!\!\begin{array}{c}
\partial w/\partial p+\partial v/\partial q\\
\partial w/\partial q-\partial v/\partial p\\
\partial w/\partial r+\partial v/\partial s\\
\partial w/\partial s-\partial v/\partial r
\end{array}\!\!\right\rgroup
\!=\!N\left(\frac{\partial}{\partial s},\frac{\partial}{\partial q}\right)$$
$$\textstyle
N\left(\frac{\partial}{\partial s},\frac{\partial}{\partial p}\right)\!=\!
\left\lgroup\!\!\begin{array}{c}
\partial v/\partial p-\partial w/\partial q\\
\partial v/\partial q+\partial w/\partial p\\
\partial v/\partial r-\partial w/\partial s\\
\partial v/\partial s+\partial w/\partial r
\end{array}\!\!\right\rgroup
\!=\!N\left(\frac{\partial}{\partial r},\frac{\partial}{\partial q}\right)
\quad
N\left(\frac{\partial}{\partial r},\frac{\partial}{\partial s}\right)\!=\!
\left\lgroup\!\!\begin{array}{c}
\displaystyle
0\\ 0\\ 0\\ 0
\end{array}\!\!\right\rgroup$$
and the proof is complete.\end{proof}

\begin{thm} If J is an integrable Hermitian structure on 
$\Omega^{\mathrm{open}}\subseteq{\mathbb{R}}^4$, then
$U=J(\partial/\partial p)$ defines a conformal foliation on
$$\Omega\cap{\mathbb{R}}^3=\{(p,q,r,s)\in\Omega\mbox{\rm\ s.t.\ }p=0\}.$$ 
\end{thm}
\begin{proof} {From} (\ref{J}) we see that $U=(u,v,w)$ and we are obliged to 
show that the equations (\ref{conformality}) hold at points where 
$(u,v,w)=(1,0,0)$. But these equations are just two of the four equations from
Lemma~\ref{nijenhuis}.
\end{proof}

Of course, the conformal foliations obtained in this way are obliged to be
real-analytic. The following theorem is a well-known result in twistor theory
(see, e.g.~\cite[Proposition~7.1.3(iii)]{thebook}). For completeness we
include a proof here based on our standard normalisation. We regard an almost
Hermitian structure $J$ on $\Omega^{\mathrm open}\subseteq S^4$ as a section of
the twistor fibration~(\ref{classical}) and write $S$ for its range.

\begin{thm}\label{twistorresult} A smooth section $J$ of\/
$\pi:{\mathbb{CP}}_3\to S^4$ defined on $\Omega$ is integrable if and only if\/
$S\equiv J(\Omega)$ is a complex submanifold of\/~${\mathbb{CP}}_3$.
\end{thm}
\begin{proof}
Without loss of generality we may suppose that the value of $J$ at some point
in $\Omega$ is given by (\ref{J}) for $(u,v,w)=(1,0,0)$ and check the statement
of the theorem above that point. The tangent space to $S$ is
given there in the local co\"ordinates (\ref{cp3coords}) by
$${\mathrm{span}}_{{\mathbb{R}}}\left\{
\left\lgroup\begin{array}c
1\\ 0\\ 0\\ 0\\ 0\\ \partial v/\partial p\\ \partial w/\partial p
\end{array}\right\rgroup,
\left\lgroup\begin{array}c
0\\ 1\\ 0\\ 0\\ 0\\ \partial v/\partial q\\ \partial w/\partial q
\end{array}\right\rgroup,
\left\lgroup\begin{array}c
0\\ 0\\ 1\\ 0\\ 0\\ \partial v/\partial r\\ \partial w/\partial r
\end{array}\right\rgroup,
\left\lgroup\begin{array}c
0\\ 0\\ 0\\ 1\\ 0\\ \partial v/\partial s\\ \partial w/\partial s
\end{array}\right\rgroup
\right\}$$
and this is preserved by ${\mathbb{J}}$ (as in (\ref{bigJ}) with
$(u,v,w)=(1,0,0)$) if and only if the equations of Lemma~\ref{nijenhuis} are
satisfied.
\end{proof}

\section{The twistor theory of conformal foliations}\label{thetwistortheory}
The following result is the direct analogue of Theorem~\ref{twistorresult} for 
conformal foliations of $\Omega^{\mathrm{open}}\subseteq S^3$ in which we 
regard a unit vector field on $\Omega$ as a section of the twistor 
fibration~(\ref{cr}).

\begin{thm}\label{mainthm}
A smooth section $U$ of $\pi:Q\to S^3$ over
$\Omega^{\mathrm{open}}\subseteq S^3$ defines a conformal foliation if and only
if its range\/ $M\equiv U(\Omega)$ is a CR submanifold of\/~$Q$ of CR 
dimension\/~$1$.
\end{thm}
\begin{proof}
Without loss of generality we may suppose that the value of $U$ at some point
in $\Omega$ is given by (\ref{thisisU}) for $(u,v,w)=(1,0,0)$ and check the 
statement of the theorem above that point. The tangent space to $M$ is
given there in the local co\"ordinates (\ref{Qcoords}) by
$${\mathrm{span}}_{{\mathbb{R}}}\left\{
\left\lgroup\begin{array}c
1\\ 0\\ 0\\ 0\\ \partial v/\partial q\\ \partial w/\partial q
\end{array}\right\rgroup,
\left\lgroup\begin{array}c
0\\ 1\\ 0\\ 0\\ \partial v/\partial r\\ \partial w/\partial r
\end{array}\right\rgroup,
\left\lgroup\begin{array}c
0\\ 0\\ 1\\ 0\\ \partial v/\partial s\\ \partial w/\partial s
\end{array}\right\rgroup
\right\}.$$
The contact form $\theta$ is simply $dq$ when $(u,v,w)=(1,0,0)$ whence $H$
intersects the tangent space to $M$ as the span of the second and third of
these three vectors and then this span is preserved by ${\mathbb{J}}$ if and
only if the equations (\ref{conformality}) characterising conformal foliations
are satisfied.
\end{proof}

\begin{cor}\label{mainresult} Let $U$ be a real-analytic unit vector field on\/
$\Omega^{\mathrm{open}}\subseteq S^3$. The following conditions are equivalent.
\begin{itemize}
\item $U$ defines a conformal foliation;
\item $U(\Omega)\subset Q$ is locally defined by a CR function;
\item $U(\Omega)=S\cap Q$ for some complex hypersurface\/
$S$ in\/ ${\mathbb{CP}}_3$;
\item $U=J(\partial/\partial p)$ for some orthogonal complex structure\/ $J$
on\/~$S^4$,
\end{itemize}
where the complex hypersurface\/ $S$ need only be defined and non-singular
near\/ $U(\Omega)\subset Q$ and, similarly, the integrable Hermitian
structure\/ $J$ need only be defined near\/ $U\subseteq S^3\subset S^4$. In the
fourth condition\/, $S^4$~is equipped with its round metric, an orientation
induced by stereographic projection, and a great subsphere\/ $S^3\subset S^4$
to which\/ $\partial/\partial p$ denotes the unit normal field.
\end{cor}
\begin{proof} In the real-analytic case, we may employ 
Theorem~\ref{holomorphicextension} to extend CR data from $Q$ to holomorphic 
data on~${\mathbb{CP}}_3$. All other equivalences have already been discussed.
\end{proof}

\section{An example---the Hopf fibration}\label{egHopf}
To employ Corollary~\ref{mainresult} in practise we need to be more specific 
about the identifications (\ref{cp3coords}) and~(\ref{Qcoords}). For 
(\ref{cp3coords}) the convenient choice is
\begin{equation}\label{convenient}
\left\lgroup\!\!\begin{array}c p+iq\\ r+is\\ u\\ v+iw\end{array}
\!\!\right\rgroup=
\frac{1}{|Z_3|^2+|Z_4|^2}\!\!\left\lgroup\!\!\begin{array}c
Z_2\bar Z_3+Z_4\bar Z_1\\
Z_1\bar Z_3-Z_4\bar Z_2\\
|Z_3|^2-|Z_4|^2\\
2iZ_4\bar Z_3
\end{array}\!\!\right\rgroup\end{equation}
for which is it evident, from (\ref{formula}) and the formula
\begin{equation}\label{stereo}{\mathbb{R}}^5=
\begin{array}c{\mathbb{R}}\\[-2pt] \oplus\\ {\mathbb{C}}^2\end{array}\supset 
S^5\ni\left\lgroup\!\!\begin{array}c t\\ \zeta_1\\ \zeta_2\end{array}\!\!
\right\rgroup\stackrel{\sigma}{\longmapsto}\frac{1}{1-t}\!\!
\left\lgroup\!\!\begin{array}c \zeta_1\\ \zeta_2\end{array}\!\!
\right\rgroup\in{\mathbb{C}}^2={\mathbb{R}}^4\end{equation}
for stereographic projection, that (\ref{cp3coords}) commutes. It is also 
routine to check that (\ref{convenient}) is holomorphic for the complex 
structure ${\mathbb{J}}$ given by (\ref{bigJ}) on ${\mathbb{R}}^4\times S^2$. 
More specifically, if we use
$$(x_1+iy_1,x_2+iy_2,x_3+iy_3)=(z_1,z_2,z_3)\mapsto[z_1,z_2,z_3,1]$$
as local co\"ordinates on ${\mathbb{CP}}_3$ then (\ref{convenient}) becomes
\begin{equation}\label{explicit}\begin{array}{rcl}
p&=&(x_2x_3+y_2y_3+x_1)/(x_3{}^2+y_3{}^2+1)\\
q&=&(x_3y_2-x_2y_3-y_1)/(x_3{}^2+y_3{}^2+1)\\
r&=&(x_1x_3+y_1y_3-x_2)/(x_3{}^2+y_3{}^2+1)\\
s&=&(x_3y_1-x_1y_3+y_2)/(x_3{}^2+y_3{}^2+1)\\
u&=&(x_3{}^2+y_3{}^2-1)/(x_3{}^2+y_3{}^2+1)\\
v&=&2y_3/(x_3{}^2+y_3{}^2+1)\\
w&=&2x_3/(x_3{}^2+y_3{}^2+1)
\end{array}\end{equation}
and it is straightforward to check that 
$${\mathbb{J}}\;\mbox{\bf Jac}=\mbox{\bf Jac}
\left[\begin{array}{cccccc}0&-1&0&0&0&0\\
1&0&0&0&0&0\\ 0&0&0&-1&0&0\\ 0&0&1&0&0&0\\ 0&0&0&0&0&-1\\ 0&0&0&0&1&0
\end{array}\right]$$
where \mbox{\bf Jac} is the Jacobian of the transformation~(\ref{explicit}). 

For (\ref{Qcoords}) it is convenient to change co\"ordinates on 
${\mathbb{CP}}_3$, writing
$$Q=\{[Z_1,Z_2,Z_3,Z_4]\in{\mathbb{CP}}_3\mbox{ s.t.\ }
Z_1\bar Z_4+Z_2\bar Z_3+Z_3\bar Z_2+Z_4\bar Z_1=0\}$$
instead of (\ref{thisisQ}), for then it is clear that $Q\setminus{\mathbb{I}}$ 
is identified as $\{p=0\}$ under~(\ref{convenient}). 

To write out the conclusions of Corollary~\ref{mainresult} more explicitly, let
us observe that (\ref{explicit}) implies that
$$\begin{array}{rcl}z_1&\!\!=\!\!&(r+is)z_3+(p-iq)\\
z_2&\!\!=\!\!&(p+iq)z_3-(r-is)\end{array}\mbox{ and }
\left\lgroup\!\!\begin{array}c u\\ v+iw\end{array}\!\!\right\rgroup=
\frac{1}{|z_3|^2+1}\left\lgroup\!\!
\begin{array}c |z_3|^2-1\\ 2i\bar z_3\end{array}\!\!\right\rgroup$$
and so if $S$ is locally written as 
$\{z\in{\mathbb{C}}^3\mbox{ s.t.\ }f(z_1,z_2,z_3)=0\}$ for some holomorphic 
function~$f$, then $U=(u,v,w)$ is defined by 
\begin{equation}\label{defofU}
\left\lgroup\!\!\begin{array}c u\\ v+iw\end{array}\!\!\right\rgroup=
\frac{1}{|z_3|^2+1}\left\lgroup\!\!
\begin{array}c |z_3|^2-1\\ 2i\bar z_3\end{array}\!\!\right\rgroup\end{equation}
where $z_3$ is the smooth function of $(q,r,s)$ implicitly defined by 
\begin{equation}\label{implicit}f((r+is)z_3-iq,iqz_3-(r-is),z_3)=0.
\end{equation}
Corollary~\ref{mainresult} says that the unit vector fields $U(q,r,s)$ obtained
in this way define real-analytic conformal foliations and that locally all
real-analytic conformal foliations arise in this way. 

As an example, if we take $f(z_1,z_2,z_3)=z_1-1$, then (\ref{implicit}) reads
$$(r+is)z_3-iq-1=0\qquad\mbox{whence}\qquad 
z_3=\frac{(1+iq)(r-is)}{r^2+s^2}$$
and so
\begin{equation}\label{clifford}
U=\left\lgroup\!\!\begin{array}c u\\ v\\ w\end{array}\!\!\right\rgroup=
\frac{1}{1+q^2+r^2+s^2}
\left\lgroup\!\!\begin{array}c 1+q^2-r^2-s^2\\ 
2(qr-s)\\ 2(qs+r)\end{array}\!\!\right\rgroup.\end{equation}
Integrating this vector field gives the mapping
$$\left\lgroup\!\!\begin{array}c q\\ r\\ s\end{array}\!\!\right\rgroup
\longmapsto\frac{1}{r^2+s^2}
\left\lgroup\!\!\begin{array}c
(1-q^2-r^2-s^2)r+2qs\\
(1-q^2-r^2-s^2)s-2qr\end{array}\!\!\right\rgroup,$$
which is the Hopf fibration $S^3\to S^2$, when viewed in stereographic
co\"ordinates. 

In this particular case there is also a global viewpoint as
follows. The zero locus of the function $f(z_1,z_2,z_3)=z_1-1$ has a 
non-singular closure in~${\mathbb{CP}}_3$, namely 
\begin{equation}\label{hopfS}
S\equiv\{[Z_1,Z_2,Z_3,Z_4]\in{\mathbb{CP}}_3\mbox{ s.t.\ }Z_1=Z_4\},
\end{equation}
which we may identify as ${\mathbb{CP}}_2$ via the mapping
$${\mathbb{CP}}_2\ni[W_1,W_2,W_3]\mapsto[W_1,W_2+W_3,W_2-W_3,W_1]
\in{\mathbb{CP}}_3.$$
If we restrict $\tau$ to $S$ viewed in this way, then (\ref{formula}) yields
$$[W]\stackrel{\tau}{\mapsto}
\frac{1}{\|W\|^2}\!\!
\left\lgroup\!\!\begin{array}c
W_2\bar W_3+W_3\bar W_2\\
|W_1|^2+|W_2|^2-|W_3|^2-W_2\bar W_3+W_3\bar W_2\\
-2W_1\bar W_3
\end{array}\!\!\right\rgroup\in
\begin{array}c{\mathbb{R}}\\[-2pt] \oplus\\ {\mathbb{C}}^2\end{array},$$
which we may, by an elementary ${\mathrm{SO}}(5)$ change of co\"ordinates,
rewrite as 
\begin{equation}\label{tau_restricted}
[W]\stackrel{\tau}{\mapsto}\frac{1}{\|W\|^2}\!\!
\left\lgroup\!\!\begin{array}c |W_1|^2+|W_2|^2-|W_3|^2\\
2W_1\bar W_3\\
2W_2\bar W_3
\end{array}\!\!\right\rgroup\in
\begin{array}c{\mathbb{R}}\\[-2pt] \oplus\\ {\mathbb{C}}^2\end{array}
={\mathbb{R}}^5.\end{equation}
{From} this viewpoint we see that there is a particular line
$${\mathbb{CP}}_1\cong\{[W_1,W_2,W_3]\in{\mathbb{CP}}_2\mbox{ s.t.\ }W_3=0\},$$
which is sent to $(1,0,0,0,0)\in S^4$, whilst on the affine co\"ordinate chart 
complementary to this line, following $\tau$ with stereographic projection 
(\ref{stereo}) gives
$$[w_1,w_2,1]\stackrel{\tau}{\mapsto}
\frac{1}{|w_1|^2+|w_2|^2+1}\!\!
\left\lgroup\!\!\begin{array}c |w_1|^2+|w_2|^2-1\\
2w_1\\ 2w_2
\end{array}\!\!\right\rgroup\stackrel{\sigma}\mapsto
\left\lgroup\!\!\begin{array}c w_1\\ w_2
\end{array}\!\!\right\rgroup.$$
In accordance with Theorem~\ref{twistorresult}, we see that the holomorphic
structure on ${\mathbb{R}}^4$ defined by (\ref{hopfS}) is just the standard
identification ${\mathbb{R}}^4={\mathbb{C}}^2$. To see the consequences of
Corollary~\ref{mainresult} in this case, note that
$$\begin{array}{rcl}S\cap Q&=&\left\{[Z]\in{\mathbb{CP}}_3\mbox{ s.t.\ }
\begin{array}{l}
Z_1\bar Z_4+Z_2\bar Z_3+Z_3\bar Z_2+Z_4\bar Z_1=0\\
Z_1=Z_4\end{array}\right\}\\
&\cong&\rule{0pt}{18pt}\big\{[W]\in{\mathbb{CP}}_2\mbox{ s.t.\ }
|W_1|^2+|W_2|^2=|W_3|^2\big\}\end{array},$$
namely the standard hyperquadric in ${\mathbb{CP}}_2$, 
whereupon (\ref{tau_restricted}) reduces to 
$$[W]\stackrel{\tau}{\mapsto}\frac{1}{2|W_3|^2}\!\!
\left\lgroup\!\!\begin{array}c
2W_1\bar W_3\\
2W_2\bar W_3
\end{array}\!\!\right\rgroup=\left\lgroup\!\!\begin{array}c
W_1/W_3\\
W_2/W_3
\end{array}\!\!\right\rgroup\in S^3\subset{\mathbb{C}}^2,$$
which induces the usual identification of this hyperquadric with the round
three-sphere. Since $\tau:S\cap Q\to S^3$ is an isomorphism in this case,
Corollary~\ref{mainresult} yields a global conformal foliation on~$S^3$. We
have already seen in local co\"ordinates (\ref{clifford}) that it is the
Clifford foliation whose integral curves define the Hopf fibration 
$S^3\to S^2$.

\section{Explicit constructions and comparisons}\label{constructandcompare}
The conclusions of Corollary~\ref{mainresult} may be written out more
explicitly as follows. In \S\ref{egHopf} we saw that the real-analytic
conformal foliations were generated from a holomorphic function
$f(z_1,z_2,z_3)$ by solving (\ref{implicit}) for $z_3$ as a smooth function of
$(q,r,s)$ and then using (\ref{defofU}) to define the unit vector
field~$U(q,r,s)$. Noting that the $z_3$-axis $\{(0,0,z_3)\}$ is the fibre of
the mapping $$\begin{array}{ccccccc}
{\mathbb{C}}^3&\hookrightarrow&{\mathbb{CP}}_3\setminus{\mathbb{I}}
&\stackrel{\tau}{\longrightarrow}&{\mathbb{C}}^2&=&{\mathbb{R}}^4\\
(z_1,z_2,z_3)&\mapsto&[z_1,z_2,z_3,1]&\mapsto&
\left\lgroup\!\!\begin{array}c
\rule{0pt}{17pt}\displaystyle\frac{z_2\bar z_3+\bar z_1}{z_3\bar z_3+1}\\
\rule{0pt}{22pt}\displaystyle\frac{z_1\bar z_3-\bar z_2}{z_3\bar z_3+1}
\end{array}\!\!\right\rgroup&=&
\left\lgroup\!\!\begin{array}c 
p+iq\\ 
\rule{0pt}{20pt}r+is\end{array}
\!\!\right\rgroup\end{array}$$
over the origin, if we seek conformal foliations near the origin in
${\mathbb{R}}^3$, then we may take the complex surface $S$ in
Corollary~\ref{mainresult} to be the graph $\{z_3=\Phi(z_1,z_2)\}$ for a
holomorphic function $\Phi=\Phi(z_1,z_2)$ defined near the origin 
in~${\mathbb{C}}^2$. In other words, we may take
$$f(z_1,z_2,z_3)=z_3-\Phi(z_1,z_2).$$
Notice from (\ref{defofU}) that precisely $(1,0,0)$ is excluded from the
possible values of $U$ at the origin (since in our formula (\ref{stereo}) for
stereographic projection, we have chosen to project from the north pole). 
Instead, we may insist that $U(0,0,0)=(-1,0,0)$ and obtain the following 
recasting of the equivalence of the first and third conditions in 
Corollary~\ref{mainresult}.
\begin{thm}\label{recast}There is a 1--1 correspondence between
\begin{itemize}
\item germs of holomorphic functions\/ $\Phi=\Phi(z_1,z_2)$ defined near and
vanishing at the origin in\/ ${\mathbb{C}}^2$; 
\item germs of real-analytic unit vector fields\/ $U=U(q,r,s)$ defined near 
and taking on the value\/ $(-1,0,0)$ at the origin in\/ ${\mathbb{R}}^3$ with 
the property that the foliation defined by\/ $U$ is conformal.
\end{itemize}
This correspondence is induced by the equations
\begin{equation}\label{defofUwithoutsubscript}
\left\lgroup\!\!\begin{array}c u\\ v+iw\end{array}\!\!\right\rgroup=
\frac{1}{|z|^2+1}\left\lgroup\!\!
\begin{array}c |z|^2-1\\ 2i\bar z\end{array}\!\!\right\rgroup\end{equation}
and
\begin{equation}\label{defofz}z=\Phi((r+is)z-iq,iqz-(r-is)),\end{equation}
where\/ $U=(u,v,w)=(u(q,r,s),v(q,r,s),w(q,r,s))$.
\end{thm}
An alternative proof in one direction may be obtained by implicit 
differentiation of the equation~(\ref{defofz}). Since $\Phi$ is holomorphic, 
the chain rule gives
$$\begin{array}{rcl}
\partial z/\partial q&=&
\big(-i+(r+is)\partial z/\partial q\big)\partial \Phi/\partial z_1
+\big(iz+iq\partial z/\partial q\big)\partial \Phi/\partial z_2\\
\partial z/\partial r&=&
\big(z+(r+is)\partial z/\partial r\big)\partial \Phi/\partial z_1
+\big(-1+iq\partial z/\partial r\big)\partial \Phi/\partial z_2\\
\partial z/\partial s&=&
\big(iz+(r+is)\partial z/\partial s\big)\partial \Phi/\partial z_1
+\big(i+iq\partial z/\partial s\big)\partial \Phi/\partial z_2
\end{array}$$
from which $\partial \Phi/\partial z_1$ and $\partial \Phi/\partial z_2$ may be
eliminated to obtain 
\begin{equation}\label{conformalPDE}
2z\frac{\partial z}{\partial q}+i(1+z^2)\frac{\partial z}{\partial r}
+(1-z^2)\frac{\partial z}{\partial s}=0\end{equation}
and it may be verified that this equation is precisely the condition that the
unit vector field $U=(u,v,w)$ defined by (\ref{defofUwithoutsubscript})
generate a conformal foliation (e.g., at $z=0$ we obtain (\ref{conformality})
as expected). It is interesting to note that the unit vector field defined by
(\ref{defofUwithoutsubscript}) is characterised, up to sign, as orthogonal to
the real and imaginary parts of the complex-valued field
$$Z\equiv
\left\lgroup\begin{array}{c}2z\\ i(1+z^2)\\ (1-z^2)\end{array}\right\rgroup.$$
Letting $\omega$ denote the equivalent $1$-form
\begin{equation}\label{omega}\omega\equiv 2z\,dq+i(1+z^2)\,dr+(1-z^2)\,ds,
\end{equation}
we find that 
\begin{equation}\label{wefindthat}\omega\wedge d\omega=2i\left(
2z\frac{\partial z}{\partial q}+i(1+z^2)\frac{\partial z}{\partial r}
+(1-z^2)\frac{\partial z}{\partial s}\right)dq\wedge dr\wedge ds
\end{equation}
and hence obtain a compact way of rewriting the equations~(\ref{conformalPDE}).
It leads to an alternative geometric interpretation of (\ref{conformalPDE}) as
follows. 
\begin{lemma}\label{frobenius} A non-zero complex-valued real-analytic
$1$-form\/ $\omega$ satisfies\/ $\omega\wedge d\omega=0$ if and only if\/
$\omega$ can be rescaled to be closed, i.e.\/ if and only if there is a smooth
function $\psi$ such that $d(e^\psi\omega)=0$.
\end{lemma}
\begin{proof} If $\omega$ is real-valued (and then need only be smooth), this
criterion is well-known. Specifically, $\omega\wedge d\omega$ is the
obstruction to Frobenius integrability of the two-dimensional distribution
defined by $\omega$ and the conclusion is evident. In the complex-valued case
it is clear from 
$$\omega\wedge d(e^\psi\omega)=\omega\wedge e^\psi(d\psi\wedge\omega+d\omega)=
e^\psi\omega\wedge d\omega$$
that the vanishing of $\omega\wedge d\omega$ is necessary. When $\omega$ is
real-analytic, the argument via Frobenius integrability applies in the
complexification (and the function $\psi$ is necessarily real-analytic back
on~${\mathbb{R}}^3$).
\end{proof}

Having chosen a real-analytic $\psi$ such that the rescaled $1$-form
$e^\psi\omega$ is closed, the Poincar\'e lemma implies that locally we may find
a real-analytic function
${\mathbb{R}}^3\supseteq{}^{\mathrm{open}}\Omega\stackrel{h}{\to}{\mathbb{C}}$
with $dh=e^\psi\omega$. In particular, since $\omega$ is null, the same is true
of $dh$, i.e.~$(dh)^2=0$. Writing $f$ and $g$ for the real and imaginary
parts of $h$, we conclude that $$\|\nabla f\|=\|\nabla g\|\quad\mbox{and}\quad
\langle\nabla f,\nabla g\rangle=0.$$
In other words, since in addition $dh$ is nowhere vanishing, we see that 
$h:\Omega\to{\mathbb{C}}$ is a {\em horizontally conformal 
mapping\/} \cite{thebook} and we have shown in Theorem~\ref{recast} and 
Lemma~\ref{frobenius} that locally all real-analytic horizontally conformal 
mappings arise in this way. 

Inspired by these formul{\ae}, we now present a direct construction of closed
null complex-valued $1$-forms on ${\mathbb{R}}^3$ avoiding the use of Frobenius
integrability. Firstly, let us observe that the implicit derivation leading to
(\ref{conformalPDE}) and (\ref{wefindthat}) may be easily accomplished as
follows. On the space ${\mathbb{C}}\times{\mathbb{R}}^3$ with co\"ordinates 
$(z,q,r,s)$ let us write
$$z_1\equiv(r+is)z-iq\quad\mbox{and}\quad z_2=iqz-(r-is)$$
and consider the $1$-form $\omega$ defined by~(\ref{omega}). One easily 
verifies that
$$\omega\wedge d\omega=2\,dz\wedge dz_1\wedge dz_2.$$
Therefore, if $z=z(q,r,s)$ is defined implicitly by 
$$z=\Phi(z_1,z_2)$$
for some holomorphic function $\Phi$ of two variables, then 
$$dz=d\Phi=
\frac{\partial\Phi}{\partial z_1}\,dz_1
+\frac{\partial\Phi}{\partial z_2}\,dz_2$$
and it follows immediately that $\omega\wedge d\omega=0$, as required.
 
Now on ${\mathbb{C}}^2\times{\mathbb{R}}^3$ with co\"ordinates $(w,z,q,r,s)$
let us write
\begin{equation}\label{z1andz2}
z_1\equiv(r+is)z-iqw\quad\mbox{and}\quad z_2\equiv iqz-(r-is)w\end{equation}
and consider the $1$-form
\begin{equation}\label{newomega}
\omega\equiv2wz\,dq+i(w^2+z^2)\,dr+(w^2-z^2)\,ds.\end{equation}
One easily verifies that 
\begin{equation}\label{domega}
d\omega=2i(dz\wedge dz_1-dw\wedge dz_2)=2i\,d(z\,dz_1-w\,dz_2).\end{equation}
Therefore, if $z=z(q,r,s)$ and $w=w(q,r,s)$ are defined implicitly by 
\begin{equation}\label{ourAnsatz}z\,dz_1-w\,dz_2=d(\Xi(z_1,z_2))=
\frac{\partial\Xi}{\partial z_1}\,dz_1+\frac{\partial\Xi}{\partial z_2}\,dz_2
\end{equation}
for some holomorphic function $\Xi$ of two variables, then $d\omega=0$ is
manifest. The formula (\ref{newomega}) for $\omega$ ensures that
$\omega^2=0$ (and also that $\omega$ is non-zero wherever $w$ is non-zero).

Postponing for the moment the precise relationship between these two 
constructions, we note that this second construction is extremely similar to 
Nurowski's method~\cite{n}, which is as follows.
With the same notation he rewrites (\ref{domega}) as
$$d\omega=2i(dz\wedge dz_1-dw\wedge dz_2)=2i\,d(z_2\,dw-z_1\,dz),$$
concluding that if $w$ and $z$ are implicitly defined by 
\begin{equation}\label{NurowskiAnsatz}
z_2\,dw-z_1\,dz=d(F(z,w))=
\frac{\partial F}{\partial z}\,dz+\frac{\partial F}{\partial w}\,dw
\end{equation}
for some holomorphic function $F$ of two variables, then~$d\omega=0$.

The relationship between these two methods is clear. In (\ref{ourAnsatz}) we 
give $(z,-w)$ in terms of $(z_1,z_2)$ whereas in (\ref{NurowskiAnsatz}) it is 
the other way round:
$$\begin{array}{rcl}
z&=&\ffrac{\partial\Xi}{\partial z_1}(z_1,z_2)\\
\tilde w&=&\rule{0pt}{20pt}\ffrac{\partial\Xi}{\partial z_2}(z_1,z_2)
\end{array}\quad\mbox{versus}\quad
\begin{array}{rclcl}
z_1&=&-\ffrac{\partial F}{\partial z}(z,w)
&=&\ffrac{\partial\tilde F}{\partial z}(z,\tilde w)\\
z_2&=&\phantom{-}\ffrac{\partial F}{\partial w}(z,w)
&=&\rule{0pt}{20pt}\ffrac{\partial\tilde F}{\partial\tilde w}(z,\tilde w),
\end{array}$$
where $\tilde w\equiv-w$ and
$\tilde F(z,\tilde w)\equiv -F(z,-\tilde w)=-F(z,w)$. In either case, the key 
point is that the Jacobian of the transformation is symmetric
\begin{equation}\label{jac}
\frac{\partial(z,\tilde w)}{\partial(z_1,z_2)}=
\left(\begin{array}{cc}
\Xi_{z_1z_1}&\Xi_{z_1z_2}\\
\Xi_{z_1z_2}&\Xi_{z_2z_2}
\end{array}\right)\quad\mbox{v}\quad
\frac{\partial(z_1,z_2)}{\partial(z,\tilde w)}=
\left(\begin{array}{cc}\tilde F_{zz}&\tilde F_{z\tilde w}\\ 
\tilde F_{z\tilde w}&\tilde F_{\tilde w\tilde w}
\end{array}\right).\end{equation}
Since the inverse of an invertible symmetric $2\times 2$ matrix is necessarily
symmetric, we may obtain one prescription from the other by inverting the
relationship between $(z,-w)$ and $(z_1,z_2)$. This assumes, of course, that
this relationship is indeed invertible: we shall come back to this point
shortly. 

Firstly, we shall explain the precise link between Theorem~\ref{recast} and the
construction determined by formul{\ae} (\ref{z1andz2}), (\ref{newomega}),
and~(\ref{ourAnsatz}). Recall that Theorem~\ref{recast} may be viewed as
generating all real-analytic null $1$-forms $\omega$ near the origin 
in~${\mathbb{R}}^3$ satisfying 
$$\omega|_{(0,0,0)}=i\,dr+ds\quad\mbox{and}\quad\omega\wedge d\omega=0.$$
Hence, to compare with (\ref{newomega}) we should suppose that $w$ and $z$ are 
given as general holomorphic functions of $(z_1,z_2)$, say
\begin{equation}\label{defoftildeS}
z=\Xi_1(z_1,z_2)\qquad w=-\Xi_2(z_1,z_2)\end{equation}
but insist that $\Xi_1(0,0)=0$ and $\Xi_2(0,0)=-1$ for then, assuming that our
construction makes sense, it will certainly create a real-analytic null
$1$-form near the origin in ${\mathbb{R}}^3$ with $\omega|_{(0,0,0)}=i\,dr+ds$
and it remains to explain the geometric origin of the equation $d\omega=0$ and
the precise link with Theorem~\ref{recast}. Regard $(z_1,z_2,z,w)$
as co\"ordinates on ${\mathbb{C}}^4$ and let
$$\eta\equiv dz\wedge dz_1- dw\wedge dz_2,$$
be a non-degenerate closed symplectic form on~${\mathbb{C}}^4\setminus\{0\}$.
Let $\tilde S$ denote the complex surface through the point
$(0,0,0,1)\in{\mathbb{C}}^4$ defined by~(\ref{defoftildeS}).
\begin{lemma}\label{Lagrange} There is a holomorphic function\/ $\Xi(z_1,z_2)$
such that locally\/ $\Xi_j=\partial\Xi/\partial z_j$ if and only
if\/~$\eta|_{\tilde S}=0$, i.e.~if and only if\/ $\tilde S$ is Lagrangian.
\end{lemma}
\begin{proof} We compute
$$\begin{array}{rcl}\eta|_{\tilde S}
&=&\left(\ffrac{\partial\Xi_1}{\partial z_1}\,dz_1
+\ffrac{\partial\Xi_1}{\partial z_2}\,dz_2\right)\wedge dz_1+
\left(\ffrac{\partial\Xi_2}{\partial z_1}\,dz_1
+\ffrac{\partial\Xi_2}{\partial z_2}\,dz_2\right)\wedge dz_2\\
&=&\rule{0pt}{20pt}\ffrac{\partial\Xi_1}{\partial z_2}\,dz_2\wedge dz_1+
\ffrac{\partial\Xi_2}{\partial z_1}\,dz_1\wedge dz_2
=\left(\ffrac{\partial\Xi_2}{\partial z_1}-\ffrac{\partial\Xi_1}{\partial z_2}
\right)dz_1\wedge dz_2,\end{array}$$
which vanishes if and only if the holomorphic $1$-form
$\Xi_1\,dz_1+\Xi_2\,dz_2$ is closed. This is the case if and only if this
$1$-form is locally of the form $d\Xi$ for some holomorphic function~$\Xi$, as
required.
\end{proof}
\noindent By construction~(\ref{defoftildeS}), the surface $\tilde S$ passes
through the point $(0,0,0,1)$ and is transverse to the $(z,w)$-plane there. It
follows that the image of a sufficiently small open subset of $\tilde S$ around
$(0,0,0,1)$ under the natural projection
${\mathbb{C}}^4\setminus\{0\}\to{\mathbb{CP}}_3$ is a complex surface 
$S\subset{\mathbb{CP}}_3$ containing the point $[0,0,0,1]$ and in the local
co\"ordinates
$${\mathbb{C}}^3\ni(z_1,z_2,z)\mapsto[z_1,z_2,z,1]\in{\mathbb{CP}}_3$$
is transverse to the $z$-axis. We may write such a surface $S$ as a
graph 
$$S=\{(z_1,z_2,z)\in{\mathbb{C}}^3\mbox{ s.t.\ }z=\Phi(z_1,z_2)\}$$
for some uniquely determined holomorphic function $\Phi(z_1,z_2)$ defined near 
and vanishing at the origin. More explicitly, there is a holomorphic function 
$\Phi(z_1,z_2)$ so that 
$$\left.\begin{array}{rcl}z&=&(\partial\Xi/\partial z_1)(z_1,z_2)\\
-w&=&(\partial\Xi/\partial z_2)(z_1,z_2)\end{array}\right\}\implies
\frac{z}{w}=\Phi\left(\frac{z_1}{w},\frac{z_2}{w}\right).$$
If we now substitute for $(z_1,z_2)$ according to (\ref{z1andz2}) and write
$z/w$ as $\zeta$, then we find that 
$$(\ref{ourAnsatz}) \implies 
\zeta=\Phi\left((r+is)\zeta-iq,iq\zeta-(r-is)\right),$$
which coincides with (\ref{defofz}) with $\zeta$ substituted for~$z$.
Moreover, we may rescale the $1$-form (\ref{newomega}) as
$$\hat\omega\equiv\frac{1}{w^2}\omega
=2\zeta\,dq+i(1+\zeta^2)\,dr+(1-\zeta^2)\,ds,$$
which coincides with (\ref{omega}) save that again $\zeta$ is substituted 
for~$z$. This is exactly as expected from Theorem~\ref{recast} and 
Lemma~\ref{frobenius} with $e^\psi=w^2$. In other words, we have found the 
precise geometric link between 
\begin{itemize}
\item the equivalence of the first and third conditions of 
Corollary~\ref{mainresult}, as recast in Theorem~\ref{recast}, and then 
rewritten in terms of a null $1$-form $\omega$ via equations~(\ref{omega})
and (\ref{wefindthat});
\item the direct construction of real-analytic closed null $1$-forms given by
formul{\ae} (\ref{z1andz2}), (\ref{newomega}), and~(\ref{ourAnsatz}).
\end{itemize}
Specifically, the complex surface $S$ appearing in the third condition of
Corollary~\ref{mainresult} is obtained as the image under the canonical
projection ${\mathbb{C}}^4\setminus\{0\}\to{\mathbb{CP}}_3$ of the Lagrangian 
surface $\tilde S$ defined by~(\ref{ourAnsatz}), namely
$$\tilde S=\left\{(z_1,z_2,z,w)\in{\mathbb{C}}^4\mbox{ s.t.\ }
z=\frac{\partial\Xi}{\partial z_1}(z_1,z_2)\enskip\mbox{and}\enskip
w=-\frac{\partial\Xi}{\partial z_2}(z_1,z_2)\right\}.$$
In fact, all real-analytic closed null $1$-forms are locally 
so obtained as follows.
\begin{thm}\label{lagrangian}There is a 1--1 correspondence between
\begin{itemize}
\item germs of complex Lagrangian submanifolds $\tilde S\subset{\mathbb{C}}^4$ 
with respect to the symplectic form\/ $\eta\equiv dz\wedge dz_1- dw\wedge dz_2$
passing through\/ $(0,0,0,1)$ and transverse to the $(z,w)$-plane there;
\item germs of real-analytic closed null\/ $1$-forms\/ $\omega$ at the origin
in\/ ${\mathbb{R}}^3$ taking on the value\/ $i\,dr+ds$ there.
\end{itemize}
This correspondence is induced by writing\/ $\tilde S$ locally as a graph
$$z=\Xi_1(z_1,z_2)\qquad w=-\Xi_2(z_1,z_2)$$
for holomorphic functions\/ $\Xi_1(z_1,z_2)$ and\/ $\Xi_2(z_1,z_2)$, using the 
equations 
$$\begin{array}{rcl}
z&=&\phantom{-}\Xi_1((r+is)z+(p-iq)w,(p+iq)z-(r-is)w)\\
w&=&-\Xi_2((r+is)z+(p-iq)w,(p+iq)z-(r-is)w)\end{array}$$
implicitly to define\/ $z=z(p,q,r,s)$ and\/ $w=w(p,q,r,s)$, restricting the
real-analytic functions\/ $z$ and\/ $w$ to\/ $\{p=0\}$, and finally setting
\begin{equation}\label{newomega-bis}
\omega=2wz\,dq+i(w^2+z^2)\,dr+(w^2-z^2)\,ds.\end{equation}
\end{thm}
\begin{proof}
Let us firstly establish a 1--1 correspondence, induced by exactly the same
procedure between
\begin{itemize}
\item germs of complex submanifolds $\tilde S\subset{\mathbb{C}}^4$ 
through $(0,0,0,1)$ and transverse to the $(z,w)$-plane there;
\item germs of real-analytic null $1$-forms $\omega$ at the origin
in ${\mathbb{R}}^3$ taking on the value $i\,dr+ds$ there and satisfying
\begin{equation}\label{heapscool}
\sigma\wedge d\omega=0\quad\forall\mbox{ $1$-forms }\sigma\mbox{ s.t.\ }
\sigma\omega=0.\end{equation}
\end{itemize}
For this, let us note that (\ref{newomega-bis}) is the general form of
a null $1$-form on ${\mathbb{R}}^3$ and that near the origin $z=z(q,r,s)$ and
$w=w(q,r,s)$ are uniquely determined by $z(0,0,0)=0$, $w(0,0,0)=1$, and
continuity. Then, the $1$-forms complex-orthogonal to $\omega$ are spanned by
$$\sigma_1\equiv z\,dq+iw\,dr+w\,ds\enskip\mbox{and}\enskip
\sigma_2\equiv w\,dq+iz\,dr-z\,ds$$
so one easily computes that (\ref{heapscool}) holds if and only if the operator
\begin{equation}\label{keyoperator}
2wz\frac{\partial}{\partial q}+i(w^2+z^2)\frac{\partial}{\partial r}
+(w^2-z^2)\frac{\partial}{\partial s}\end{equation}
annihilates both $z(q,r,s)$ and $w(q,r,s)$. Now let us ask what it means for a 
surface $\tilde S$ through $(0,0,0,1)$ in ${\mathbb{C}}^4$ to be complex in 
terms of the local co\"ordinates $(\alpha,\beta,z,w)$ defined in terms of 
$(z_1,z_2,z,w)$ by the relations
$$z_1=\alpha z+\bar\beta w\quad\mbox{and}\quad z_2=\beta z-\bar\alpha w,$$
which are obtained by substituting
$$\alpha\equiv r+is\quad\beta\equiv p+iq\quad z_1\equiv Z_1\quad z_2\equiv Z_2
\quad z\equiv Z_3\quad w\equiv Z_4$$
into (\ref{convenient}). Certainly, we may write $\tilde S$ locally as a smooth
graph
$$z=z(\alpha,\beta)\qquad w=w(\alpha,\beta).$$
We may then verify, using the chain rule to change co\"ordinates, that $\tilde
S$ is complex if and only if $z(\alpha,\beta)$ and $w(\alpha,\beta)$ are
annihilated by the operators
\begin{equation}\label{thefullCauchyRiemannequations}
w\frac{\partial}{\partial\alpha}-z\frac{\partial}{\partial\bar\beta}
\quad\mbox{and}\quad
w\frac{\partial}{\partial\beta}+z\frac{\partial}{\partial\bar\alpha}.
\end{equation}
Recall that, with the conventions of~\S\ref{egHopf}, the hyperquadric 
$Q\subset{\mathbb{CP}}_3$ is covered by 
$$\tilde Q\equiv\{(z_1,z_2,z,w)\in{\mathbb{C}}^4\mbox{ s.t.\ }
z_1\bar w+z_2\bar z+z\bar z_2+w\bar z_1=0\},$$
a CR hypersurface in ${\mathbb{C}}^4\setminus\{0\}$ of Levi-signature
$(+,0,-)$. Sitting over Theorem~\ref{holomorphicextension} and similarly
proved, suppose $\tilde M\subset\tilde\Omega^{\mathrm{open}}\subseteq\tilde 
Q\subset{\mathbb{C}}^4$ is a real-analytic submanifold of real dimension $3$ 
and CR dimension~$1$. Then $\tilde M$ extends into ${\mathbb{C}}^4$ as a 
complex submanifold $\tilde S$ and this extension is germ-unique. Since,
$$z_1\bar w+z_2\bar z+z\bar z_2+w\bar z_1=2(|z|^2+|w|^2)p$$
we conclude that $\tilde Q$ is defined as the zero locus of $p$ in a
neighbourhood of~$(0,0,0,1)$. Writing (\ref{thefullCauchyRiemannequations})
more fully, we see that $\tilde S$ is complex if and only if the operators
$$w\frac{\partial}{\partial r}-iw\frac{\partial}{\partial s}
-z\frac{\partial}{\partial p}-iz\frac{\partial}{\partial q}
\quad\mbox{and}\quad
w\frac{\partial}{\partial p}-iw\frac{\partial}{\partial q}
+z\frac{\partial}{\partial r}+iz\frac{\partial}{\partial s}$$
annihilate both $z(p,q,r,s)$ and $w(p,q,r,s)$. On $\tilde Q=\{p=0\}$ it follows
that the operator (\ref{keyoperator}) annihilates both $z(0,q,r,s)$ and
$w(0,q,r,s)$. We have shown that the complex submanifold $\tilde S$ gives rise
to a real-analytic null $1$-form $\omega$ satisfying~(\ref{heapscool}).
Conversely, it is readily verified that being annihilated by
(\ref{keyoperator}) is exactly the condition that the functions $z(q,r,s)$ and
$w(q,r,s)$ define a CR submanifold $\tilde M$ of~$\tilde Q$. (Warning: this is
not to say that the defining functions 
$$(q,r,s,z,w)\mapsto z-z(q,r,s)\quad\mbox{and}\quad
(q,r,s,z,w)\mapsto w-w(q,r,s)$$
are CR functions on $\tilde Q$ but only that the tangent space defined by 
them intersects the contact distribution defined by
$$(|w|^2-|z|^2)\,dq+i(z\bar w-w\bar z)\,dr-(z\bar w+w\bar z)\,ds$$
in a complex subspace.) When $z(q,r,s)$ and $z(q,r,s)$ are real-analytic this
CR submanifold $\tilde M$ extends germ-uniquely into ${\mathbb{C}}^4$
as~$\tilde S$, a complex surface: we have now shown the equivalence of the two
entities claimed to be equivalent at the beginning of this proof.

To finish the proof it remains to show that $d\omega=0$ if and only if
$\tilde S$ is Lagrangian. Since $\eta$ is a holomorphic form of type $(2,0)$
its pullback $\eta|_{\tilde S}$ is a holomorphic section of the canonical
bundle of~$\tilde S$. Hence $\eta|_{\tilde S}$ vanishes near $\tilde M$ if and 
only if $\eta|_{\tilde M}=0$. However, as we have already implicitly noticed 
in~(\ref{domega}), writing $\tau:{\tilde M}\to{\mathbb{R}}^3$ for the 
canonical projection,
$$\tau^*(d\omega)=2i\eta|_{\tilde M}$$
whence $\eta|_{\tilde M}=0$ if and only if $\omega$ is closed.
\end{proof}

As an example of Theorem~\ref{lagrangian} in action, let us consider the
$1$-form
\begin{equation}\label{eg}\omega\equiv(1+ir+s)(i\,dr+ds)\end{equation}
on ${\mathbb{R}}^3$. Evidently, it is closed, null, and real-analytic. Near 
the origin, it is of
the form (\ref{newomega-bis}) for $z(q,r,s)\equiv 0$ and
$w(q,r,s)=\sqrt{1+ir+s}$ where $\sqrt{\phantom{x}}$ is a branch of square root 
with $\sqrt{1}=1$. Theorem~\ref{lagrangian} says firstly that there are 
holomorphic functions $\Xi_1(z_1,z_2)$ and $\Xi_2(z_1,z_2)$ such that 
$$\begin{array}{rcl}
0&=&\phantom{-}\Xi_1((p-iq)\sqrt{1+ir+s},-(r-is)\sqrt{1+ir+s})\\
\sqrt{1+ir+s}&=&-\Xi_2((p-iq)\sqrt{1+ir+s},-(r-is)\sqrt{1+ir+s})\end{array}$$
and this is clear by taking $\Xi_1\equiv 0$ and $\Xi_2(z_1,z_2)=-g(z_2)$ where 
$g(z_2)$ is implicitly defined near $z_2=0$ by
$$\zeta=g(i(\zeta^2-1)\zeta)\enskip\mbox{near }
\zeta=1.$$
Additionally, in conjunction with Lemma~\ref{Lagrange},
Theorem~\ref{lagrangian} says that we can find a function $\Xi(z_1,z_2)$ such
that $\Xi_j=\partial\Xi/\partial z_j$ and in this example we may take
$\Xi(z_1,z_2)=-f(z_2)$ where $f$ is any primitive for~$g$.

Returning now to the relationship between our construction, now formulated in
Theorem~\ref{lagrangian}, and Nurowski's construction in~\cite{n}, we see that
\cite{n} provides an alternative way of writing a real-analytic closed null
$1$-form $\omega$ near the origin if and only if the left hand matrix from
(\ref{jac}) is invertible. {From} (\ref{newomega}) and~(\ref{ourAnsatz}),
however, one easily computes that
$$dz\wedge dw=\omega\wedge (\Xi_{11}\Xi_{22}-\Xi_{12}{}^2)dq\quad\mbox{at the 
origin}.$$
It follows that \cite{n} pertains if and only if $dz\wedge dw$ is non-vanishing
at the origin (a requirement independent of choice of co\"ordinates on
${\mathbb{R}}^3$). Generically, this is true but not so for (\ref{eg}) nor for 
$$\omega=\frac{i\,dr+ds}{(1+(ir+s))^2},$$
(which gives the Hopf fibration in another guise). To repair Nurowski's
construction one can treat the case $dz\wedge dw=0$ separately or one can view
his construction as giving all real-analytic null $1$-forms in a neighbourhood
of $\infty\in S^3$.

In \cite[Example~2.6]{bp} it is explained how to use special holomorphic
co\"ordinates on the surface $S\hookrightarrow{\mathbb{CP}}_3$ from
Corollary~\ref{mainresult} to produce the general real-analytic horizontally
conformal submersion on~${\mathbb{R}}^3$. These special co\"ordinates are
adapted to the contact structure on ${\mathbb{CP}}_3$ induced by our symplectic
form~$\eta$. It is unclear whether one can use the lifted Lagrangian surface
$\tilde S\hookrightarrow{\mathbb{C}}^4$ from Theorem~\ref{lagrangian} to
generate such adapted co\"ordinates directly.

\section{A counterexample---the eikonal equation}
At the end of \S\ref{CRgeometry} we mentioned that there are smooth
$3$-dimensional CR submanifolds $M\subset Q$ of CR dimension 1 that are not
real-analytic. In view of Theorem~\ref{mainthm}, to find such an $M$ it
suffices to find a smooth conformal foliation of
$\Omega^{\mathrm{open}}\subseteq{\mathbb{R}}^3$ that is not real-analytic. To
construct such a foliation, let $\phi:{\mathbb{R}}^2\to{\mathbb{R}}$ be any 
smooth function, let 
$$\Gamma\equiv\{(r,s)\in{\mathbb{R}}^2\mbox{ s.t.\ }s=\phi(r)\}$$
denote its graph, and define 
$$\rho(r,s)=\Big\{\!\!\begin{array}l
\mbox{$\phantom{{}-{}}$distance from $(r,s)$ to~$\Gamma$, if $s\geq\phi(r)$}\\
\mbox{${}-{}$distance from $(r,s)$ to~$\Gamma$, if $s\leq\phi(r)$}.
\end{array}$$
Evidently, the function $\rho$ is a smooth solution of the eikonal equation
$\|\nabla\rho\|=1$ in a suitable neighbourhood ${\mathcal{N}}$ of~$\Gamma$. 
It follows that 
$$\Omega={\mathbb{R}}\times{\mathcal{N}}\stackrel{\pi}{\longrightarrow}
{\mathbb{R}}^2\quad\mbox{given by}\quad\pi(q,r,s)=(q,\rho(r,s))$$
is a horizontally conformal submersion. The corresponding conformal foliation
is real-analytic if and only if the same is true for our original
function~$\phi$.

\section{Further equivalences---the Kerr Theorem}\label{kerr}
In this section we relate the notions of conformal foliation and shear-free ray
congruence from relativity. In \cite{bw} this relationship was used to derive
the twistor description of conformal foliations,
i.e.~Corollary~\ref{mainresult}.

We may regard Euclidean space ${\mathbb{R}}^3$ as the slice $\{t=0\}$ in
Minkowski space ${\mathbb{R}}^{3,1}$ equipped with co\"ordinates $(q,r,s,t)$
and pseudo-metric
$$dq^2+dr^2+ds^2-dt^2$$
in the usual way. If $U$ is a smooth unit vector field on
$\Omega^{\mathrm{open}}\subseteq{\mathbb{R}}^3$, then $U+\partial/\partial t$
is a null direction in ${\mathbb{R}}^{3,1}$ defined and smoothly varying along
$U\subset{\mathbb{R}}^{3,1}$. We consider the region swept out in
${\mathbb{R}}^{3,1}$ by the null rays emanating from $\Omega$ in the direction
given by $U+\partial/\partial t$. For a suitable neighbourhood $\tilde\Omega$
of $\Omega$ in ${\mathbb{R}}^{3,1}$, this construction gives what is called a
`ray congruence' in the relativity literature. It is a smooth family of null
geodesics with one such geodesic passing through each point. Locally all such
congruences near the slice $\{t=0\}$ arise in this way and the vector field
$U$ defines a conformal foliation if and only if the corresponding ray
congruence is `shear-free'~\cite{NT} (the `shear' being a measure of the
distortion of circles to ellipses in the normal bundle to the foliation of
$\Omega$ defined by~$U$).

From Corollary~\ref{mainresult} we conclude that the real-analytic shear-free
ray congruences defined near $\{t=0\}$ are locally in $1$--$1$ correspondence
with complex hypersurfaces $S\subset{\mathbb{CP}}_3$ meeting
$Q\subset{\mathbb{CP}}_3$ as discussed in \S\ref{thetwistortheory}. This is the
Kerr Theorem~\cite[Theorem~7.4.8]{NT}. Furthermore, the smooth shear-free ray
congruences correspond to CR submanifolds of $Q$ of CR dimension~$1$ as in
Theorem~\ref{mainthm}. This fact is also observed on pages 220--222
of~\cite{NT}. As detailed in~\cite{NT}, the Kerr Theorem was highly 
instrumental in Penrose's development of twistor theory.


\begin{thebibliography}{X}
    
\bibitem{a} M.F. Atiyah,
{\em Geometry of Yang--Mills Fields},
Scuola Normale Superiore Pisa 1979. 

\bibitem{ber} M.S. Baouendi, P. Ebenfelt, and L.P. Rothschild,
{\em Real Submanifolds in Complex Space and their Mappings},
Princeton University Press 1999. 

\bibitem{bp} P. Baird and R. Pantilie,
{\em Harmonic morphisms on heaven spaces},
Bull. Lond. Math. Soc. {\bf 41} (2009) 198-Ð204. 

\bibitem{bw} P. Baird and J.C. Wood,  
{\em Harmonic morphisms and shear-free ray congruences} (2002), 
\verb+www.maths.leeds.ac.uk/pure/staff/wood/BWBook/BWBook.html+

\bibitem{thebook} P. Baird and J.C. Wood,
{\em Harmonic Morphisms between Riemannian Manifolds},
Oxford University Press 2003.

\bibitem{h} L. H\"ormander, 
{\em An Introduction to Complex Analysis in Several Variables},
Van Nostrand 1966, North-Holland 1973. 

\bibitem{n} P. Nurowski,
{\em Construction of conjugate functions},
Ann. Glob. Anal. Geom. {\bf 37} (2010) 321--326.

\bibitem{NT} R. Penrose and W. Rindler,
{\em Spinors and Space-time vol.~2},
Cambridge University Press 1986.

\end{thebibliography}
\end{document}